\documentclass{amsart}
\usepackage{amsfonts}
\usepackage[section]{placeins}
\usepackage{comment}
\usepackage{amsmath}
\usepackage{amssymb}
\usepackage{amsthm}
\usepackage{array}
\usepackage{graphicx} 
\usepackage{pstricks}
\usepackage{hyperref}
\usepackage{pifont}
\usepackage{xcolor}
\newcommand{\nomark}{\colorbox{red}{\textcolor{white}{\ding{55}}}}
\newcommand{\yesmark}{\colorbox{green!60!black}{\textcolor{white}{\ding{51}}}}

\def\Erdos{Erd\H{o}s}

\def\Z{{\mathbb Z}}
\def\R{{\mathbb R}}

\def\p{{\mathfrak P}}

\def\onehalf{\frac 1 2} 
\def\mod{\mbox{\ mod\ }}

\newtheorem{theorem}{Theorem}
\newtheorem{lemma}{Lemma}
\newtheorem{corollary}{Corollary}

\theoremstyle{definition} 
\newtheorem{definition}{Definition}
\newtheorem*{remark}{Remark}

\def\Erdos{Erd\H{o}s}

\def\Z{{\mathbb Z}}

\def\p{{\mathfrak P}}

\def\onehalf{\frac 1 2} 
\def\mod{\mbox{\ mod\ }}

\begin{document}
\title{Tiling an Equilateral Triangle}
\author{Michael Beeson}
\date{\today}
 \maketitle
 
\begin{abstract}
Let $ABC$ be an equilateral triangle.  For certain 
triangles $T$ (the ``tile'') and certain $N$, it is 
possible to cut $ABC$ into $N$ copies of $T$.  It is 
known that only certain shapes of $T$ are possible, but
until now very little was known about the possible values
of $N$.  Here we  prove that for $N>3$,
 $N$ cannot be prime, and study more closely the possible
 tilings when the tile has a $\pi/3$ angle.   
 All theorems have been formalized in Lean.

\end{abstract}

\section{Introduction}
The subject of this paper is $N$-tilings of the equilateral 
triangle.  More generally, triangle $ABC$ is said to be $N$-tiled
by a triangle (the ``tile'') with 
angles $(\alpha,\beta,\gamma)$,  if $ABC$ can be cut into $N$
smaller triangles congruent to the tile.  In this paper, we
restrict attention to the case of $ABC$ equilateral. 
A few pictures of $N$-tilings of equilateral 
triangles are given in the figures.

\begin{figure}[ht]   
\caption{A 3-tiling, a 6-tiling, and a 16-tiling}
\label{figure:nsquared}

\begin{center}
\includegraphics[width=0.3\textwidth]{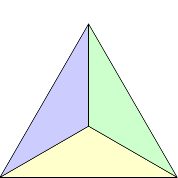}
\includegraphics[width=0.3\textwidth]{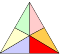}
\includegraphics[width=0.3\textwidth]{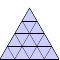}
\end{center}
\end{figure}

\begin{figure}[ht]    
\caption{A 27-tiling due to Major MacMahon 1921, rediscovered 2011}
\label{figure:prime27}
\begin{center}
\includegraphics[width=0.5\textwidth]{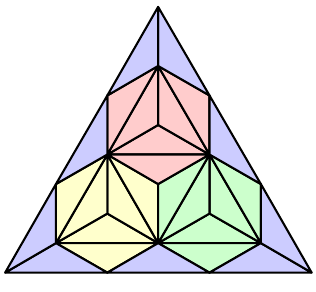}
\end{center}
\end{figure}

\begin{figure}[ht]    
\caption{$3m^2$ (hexagonal) tilings for $m=4$ and $m=5$}
\label{figure:hexagonaltilings}
\begin{center}
\includegraphics[width=0.45\textwidth]{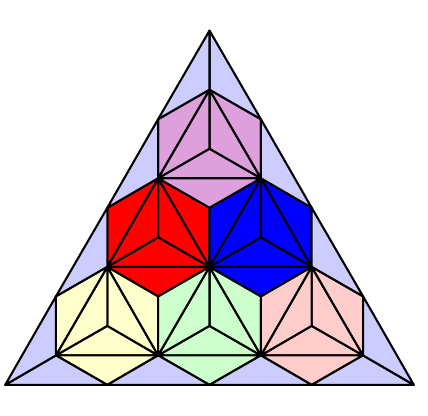}
\includegraphics[width=0.45\textwidth]{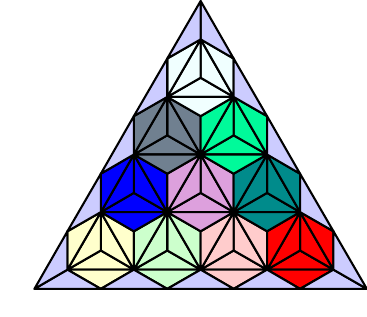}
\end{center}
\end{figure}

\begin{figure}[ht]
\caption{$N=10935$.  The tile is $(3,5,7)$ and $\gamma = 2\pi/3$. }
\label{figure:10935}
\begin{center}
\psset{unit=1cm}
\includegraphics[width=\textwidth]{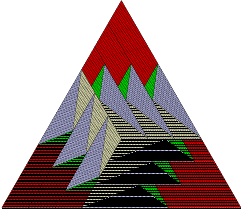}  
\end{center}
\end{figure}

\begin{figure}[ht]
\caption{$N=1215$.  The tile is $(3,5,7)$ and $\gamma = 2\pi/3$. }
\label{figure:1215}
\begin{center}
\psset{unit=1cm}
\includegraphics[width=\textwidth]{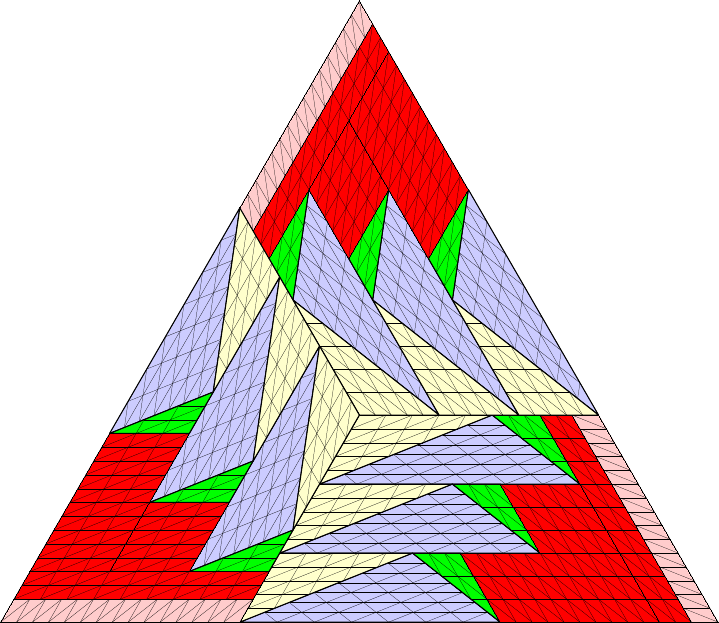} 
\end{center}
\end{figure}

Mikl\'os Laczkovich wrote many papers on the subject of tilings, 
considering not only triangles but convex polygons, and considering
not only tilings by congruent tiles, but by similar tiles as well. 
In \cite{laczkovich1995} and \cite{laczkovich2012} (especially Theorem~3.3 of \cite{laczkovich2012})
he narrowed the list of possible tiles that can be used to tile 
some equilateral triangle.  He proved that the tile $(\alpha,\beta,\gamma)$
(in some order) must be one of the following
\medskip

(i) $(\pi/3,\pi/3,\pi/3)$  (equilateral)

(ii) $(\pi/6,\pi/6,2\pi/3)$
\smallskip

(iii) $(\pi/6, \pi/2, \pi/3)$
\smallskip

(iv) $(\alpha, \beta, \pi/3)$ with $\alpha$ not a rational multiple of $\pi$
\smallskip

(v)  $(\alpha,\beta,2\pi/3)$ with $\alpha$ not a rational multiple of $\pi$
\smallskip

The possible $N$ in the first three cases are known.  In case (i),
$N$ must be a square, and any square corresponds to a tiling.  (That 
is true for any triangle $ABC$ and a tile similar to $ABC$.) 
For a proof see \cite{snover1991} or \cite{soifer}.  In cases (ii) and (iii),
$N$ must have the form $3n^2$ or $6n^2$ respectively \cite{beeson-noseven}. 
None of these can be prime except for case (ii), when $N=3$ is possible.
That leaves the last two cases as the focus of this paper. 

It may come as a surprise to the reader that tilings falling under
cases (iv) and (v) do exist.  
In 1995, Laczkovich already gave a method for constructing a tiling
of {\em some} equilateral triangle from any tile satisfying (iv) or (v),
but the $N$ involved can be quite large.   Following Laczkovich's method,
in 2018 I constructed a tiling with
 $\gamma = 2\pi/3$ and 
$N= 10935$; see Fig.~\ref{figure:10935}.  In January 2024, Bryce Herdt constructed
the 1215-tiling shown in Fig.~\ref{figure:1215}.  He retiled the larger tiling
with triangles three times longer, which is possible if one splits off a 
certain parallelogram (pink in the figure) and tiles it with a 
different orientation.
In \cite{beeson-noseven}, it is shown that $N$ must be at least 12.
There is still a gap between 12 and 1215: are there equilateral tilings with $\gamma = 2\pi/3$
and $N$ in this range?   

For tilings when the 
tile has a $\pi/3$ angle,  the smallest one I could 
calculate following Laczkovich has more than five million tiles.
But in 2024, Bryce Herdt produced better tilings.  He constructed
a 1944-tiling by $(3, 8, 7)$ of an equilateral triangle of side 216,
and a 1440-tiling by $(5,8,7)$ of an equilateral triangle with side 240.
Pictures are given in an Appendix.
 
In this paper, we prove that if an equilateral triangle is $N$-tiled according
to cases (iv) or (v) above, 
then $N$ is not a prime.  
Laczkovich  proved in  Theorem~3.3 of \cite{laczkovich2012}
that in those two cases,  the tile is rational; that is, the 
ratios of the sides of the tile are rational.  This gives us a 
good starting point.  On that basis, we proceed in case (iv) by 
elementary, if somewhat intricate, algebraic computations, arriving
at some equations that must be satisfied by $N$.  While we could not 
reduce these equations to an elegant number-theoretical criterion, 
at least we can show that $N$ cannot be prime. The algebraic 
calculations have been performed or 
checked using the computer algebra system SageMath \cite{Sage},
but they are presented here in human-checkable detail.   
\FloatBarrier

To deal with case (v), when the tile has a $2\pi/3$ angle,
we make use of a beautiful method introduced by Laczkovich in 
his 2012 paper \cite{laczkovich2012}.  In this paper, we 
give the required definitions and state the two lemmas we need,
but for the proofs of those lemmas, we refer to the cited paper. 

The result of this paper,  that equilateral triangles cannot be 
$N$-tiled when $N$ is a prime larger than 3, answers a question 
posed in \cite{soifer}.   Namely,  until now it was not known whether
there are arbitrarily large $N$ such that no equilateral triangle can 
be $N$-tiled.   Now, in view of Euclid's theorem that there are infinitely
many primes, we know the answer. 

\section{The coloring equation}
\begin{figure} [ht]
\caption{A tiling colored so that touching tiles have different colors.}
\label{figure:coloring}
\begin{center} 
\includegraphics[width=0.5\textwidth]{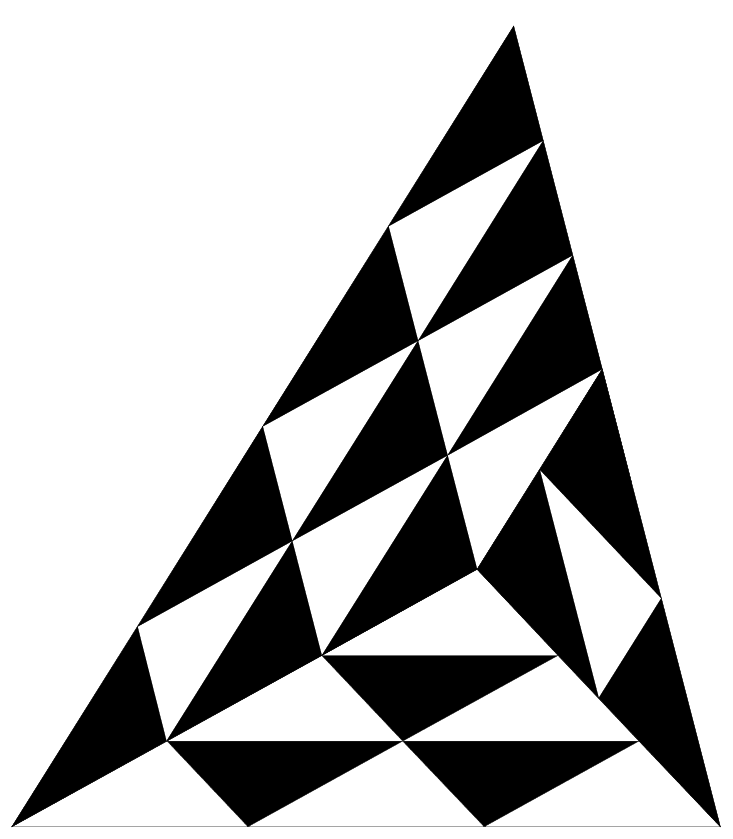}
\end{center}
\end{figure}

In this section we introduce a tool that is useful for some, but not all, 
tiling problems.  Suppose that triangle $ABC$ is tiled
by a tile with 
angles $(\alpha,\beta,\gamma)$ and sides $(a,b,c)$, and fix one
particular tile, say the one at vertex $A$ with an edge on $AB$.
 We color that tile black, 
and then we color each tile black or white, changing colors as 
we cross tile boundaries.  Under certain conditions this coloring
can be defined unambiguously, and then, we define the ``coloring number''
to be the number of black tiles minus the number of white tiles.
An example
of such a coloring is given in Fig.~\ref{figure:coloring}.  The
tiling in that figure is not a tiling of an equilateral triangle
of the type mentioned in the next theorem, because we do not 
have small examples of such tilings, and the possibility of 
coloring applies to other shapes as well.

\begin{theorem} \label{theorem:coloring}
Let the equilateral triangle $ABC$ be $N$-tiled by a tile with angles
$(\alpha,\beta,\pi/3)$ and sides $(a,b,c)$, where $\alpha$ is not a rational
multiple of $\pi$, and let $X$ be the side of $ABC$. Then the tiles can be
colored black and white so that any two tiles whose boundaries share a segment
of positive length have different colors, and the tile at $A$ is black.
For any such coloring, let $M$ be the number of black tiles minus the number
of white tiles. Then
$$ M(a+b+c) = 3X, $$
and the coloring is uniquely determined once the color of one tile is fixed.
\end{theorem}

\begin{remark}
More generally, suppose a triangle $ABC$ with sides $X$, $Y$, $Z$, where $Y$ is
opposite $A$, is tiled by the tile $(a,b,c)$ so that (i) there is just one tile
at $A$; (ii) an odd number of tiles meet at every boundary vertex; (iii) an even
number of tiles meet at every interior vertex; and (iv) the numbers of tiles at
$B$ and at $C$ are both even or both odd. Then the same kind of coloring exists,
and $X \pm Y + Z = M(a+b+c)$, with the sign $+$ or $-$ according as the number
of tiles at $B$ and $C$ is odd or even. Here a ``boundary vertex'' lies on the
boundary of $ABC$ or on an edge of another tile, and an ``interior vertex'' is
one where the angles sum to $2\pi$. We use only the case stated in
Theorem~\ref{theorem:coloring}, and only that case is formalized.  
In the theorem,  $X \pm Y + Z$ becomes $X+Y+Z$.  The reason for formalizing
only one case is to avoid needing topology to prove that the
 color so defined 
is independent of the path chosen to reach the tile from $A$.  
Informally, we prove that the parity of the number of lines crossed does
not change under a homotopy of the path, because of the hypothesis that
an even number of tiles meet at interior vertices.  Formally, that is 
slightly tricky as the homotopy passes through a vertex where the crossing 
number is not even defined.  We therefore give a more algebraic proof,
which only works for equilateral tilings.
\end{remark}

\begin{proof}[Proof of Theorem~\ref{theorem:coloring}]
We prove only the equilateral case.
  Each tile is colored black or white
according as the number of tile boundaries crossed in reaching it 
from $A$ without passing through a vertex is even or odd.  

To see that the color does not depend on the path chosen,  orient the boundary
of each tile counterclockwise, and let
$G = \{ j\pi/3 + k\alpha : j, k \in \Z \}$, taken mod $2\pi$.
Since $\alpha$ is not a rational multiple of $\pi$, each $\theta \in G$
has a unique representation $\theta = j\pi/3 + k\alpha$ with $j$ taken
mod 6; we call $j+k \pmod 2$ the {\em parity} of $\theta$.
Going around a tile counterclockwise, the direction of the edges turns
successively by the exterior angles $\pi - \alpha$, $\pi - \beta = \pi/3 + \alpha$,
and $2\pi/3$, and each of these is in $G$ with even parity. Hence if one
edge of a tile has direction in $G$, so do the other two, and all three
have the same parity; we call it the parity of the tile.
Two tiles whose boundaries share a segment of positive length traverse it
in opposite directions, so their directions there differ by $\pi = 3\cdot\pi/3$,
which has odd parity. The tile at $A$ has an edge on $AB$, of direction $0$.
Therefore, along any path from $A$, every tile reached has its edge
directions in $G$, and each boundary crossed changes the parity of the tile.
Hence the number of boundaries crossed in reaching a tile $T$ is even or odd
according as $T$ has the same parity as the tile at $A$ or not; this does not
depend on the path. Every tile can be reached from $A$ by such a path, 
since $ABC$ with the finitely many vertices removed is still connected.
The coloring so defined
gives different colors to tiles sharing a segment. Conversely, any coloring
with that property and with the tile at $A$ black changes color at each
crossing, so it is this coloring; hence $M$ is well defined.

The total length of black edges, minus the total length of white 
edges,  is $M(a+b+c)$, since $a+b+c$ is the perimeter of each tile.
Each point of an interior segment, other than a vertex, lies on the boundary 
of exactly two tiles, which share a segment there and so have different colors; 
hence the black and white lengths along each interior segment are equal.
Therefore only the edges on the boundary of $ABC$ contribute. 
A tile with an edge on $AB$ has that edge in direction $0$, so it has the
parity of the tile at $A$ and is black. Similarly the edges on $BC$ and $CA$
have directions $2\pi/3$ and $4\pi/3$, both of even parity, so every tile
with an edge on the boundary of $ABC$ is black. 
Hence the difference in the total length of black and white edges
is $X + Y + Z$. 
\end{proof}
 
\section{A tile with an angle $\gamma = \pi/3$ and $\alpha/\pi$ irrational}
In this case we have $\alpha + \beta = 2\pi/3$.  Then the possible ways
to write $\pi$ and $2\pi$ as an integer linear combination of $(\alpha,\beta,\gamma)$
are these:  
\begin{eqnarray*}
\pi &=& \alpha+\beta + \gamma \\
2\pi &=& 2\alpha + 2\beta+ 2\gamma\\
2\pi &=& 6 \gamma \\
2 \pi &=& \alpha + \beta + 4 \gamma \\
2\pi &=& 2 \alpha + 2\beta + 2 \gamma
\end{eqnarray*}
 Hence every vertex with total angle $\pi$ has an odd number
of tiles sharing that vertex  and every vertex
with total angle $2\pi$ has an even number of tiles 
sharing that vertex.  Moreover at each vertex of $ABC$ there is just one tile, with 
its $\gamma$ angle at the vertex.  Thus the coloring 
theorem, Theorem~\ref{theorem:coloring}, applies. 
It tells that,  
 when the tiles are colored black and white with the 
vertex at $A$ black, if $M$ is the number of black tiles minus the 
number of white tiles, 
\begin{eqnarray}
 M(a+b+c) &=& 3X \label{eq:coloring}
\end{eqnarray}
where $X$ is the length of each side of $ABC$.

Our second tool is the area equation, obtained by equating 
the area of $ABC$ to $N$ times the area of the tile.  Twice the 
area of $ABC$ is $X^2 \sin \gamma$, and twice the area of the tile 
is $ab \sin \gamma$, so the area equation is 
\begin{eqnarray}
X^2 &=& Nab  \label{eq:area}
\end{eqnarray}

These tools enable us to formulate a necessary condition for such 
a tiling to exist, and characterize the tile, i.e., compute 
the tile from $N$ and the coloring number of the tiling. 

\begin{theorem}\label{theorem:equilateral3}
Let triangle $ABC$ be equilateral.  Suppose it is $N$-tiled using
a tile with angles $(\alpha,\beta,\frac \pi 3)$, where $\beta$
is not a rational multiple of $\pi$.  Let $M$ be the
coloring number of the tiling.  Then
\smallskip

(i) $N$ and $M$ determine the shape of the tile uniquely; that is, 
there are algebraic formulas for $(a/c)$ and $(b/c)$ in terms of $N$ and $M$.
Specifically, $a/c$ and $b/c$ are given by
\begin{eqnarray*}
  \frac 1 2 \left( \frac {3N+M^2}{3N-M^2}\pm \frac {\sqrt{(9N-M^2)(N-M^2)}}
  {3N-M^2}\right)
\end{eqnarray*}
\smallskip

(ii)  $M^2 < N$
\end{theorem}

\begin{proof}
Let $(a,b,c)$ be the sides of the tile, opposite the angles
$(\alpha,\beta,\pi/3)$, and let $X$ be the side of $ABC$.
Since $\cos(\pi/3) = 1/2$, the law of cosines gives
\begin{eqnarray}
c^2 &=& a^2 + b^2 - ab \nonumber\\
c^2 &=& (a+b)^2 - 3ab \nonumber\\
3ab &=& (a+b)^2 - c^2 \label{eq:3ab}
\end{eqnarray}
In particular $a + b > c$.  From the coloring and area equations we have
\begin{eqnarray*}
M^2(a+b+c)^2 &=& 9X^2  \mbox{\hskip3cm squaring (\ref{eq:coloring})}\\
&=& 9Nab \mbox{\hskip2.8cm by (\ref{eq:area})}\\
&=& 3N\big((a+b)^2 - c^2\big) \mbox{\hskip0.8cm \ by (\ref{eq:3ab})}\\
&=& 3N(a+b+c)(a+b-c) \\
M^2(a+b+c) &=& 3N(a+b-c) \mbox{\hskip1.5cm  dividing by $a+b+c > 0$}\\
(3N-M^2)(a+b) &=& (3N+M^2)\,c \mbox{\hskip1.7cm  collecting terms}
\end{eqnarray*}
The right side is positive, so $3N - M^2 > 0$.  Let $s = (a+b)/c$.  Then
\begin{eqnarray}
s &=& \frac{3N+M^2}{3N-M^2} \label{eq:s}
\end{eqnarray}
Dividing (\ref{eq:3ab}) by $3c^2$, we have
\begin{eqnarray}
\frac{ab}{c^2} &=& \frac{s^2-1}{3} \nonumber\\
&=& \frac{(3N+M^2)^2-(3N-M^2)^2}{3(3N-M^2)^2}
   \mbox{\qquad by (\ref{eq:s})} \nonumber\\
&=& \frac{4M^2N}{(3N-M^2)^2} \label{eq:abc2}
\end{eqnarray}
Therefore
\begin{eqnarray*}
\left(\frac{a-b}{c}\right)^2 &=& s^2 - \frac{4ab}{c^2} \\
&=& \frac{(3N+M^2)^2 - 16M^2N}{(3N-M^2)^2}
   \mbox{\qquad by (\ref{eq:s}) and (\ref{eq:abc2})}\\
&=& \frac{9N^2 - 10M^2N + M^4}{(3N-M^2)^2} \\
&=& \frac{(9N-M^2)(N-M^2)}{(3N-M^2)^2}
\end{eqnarray*}

{\em Part (ii).}  We have $a \ne b$, since otherwise $\alpha = \beta$
by the law of sines, so $\alpha = \pi/3$, contrary to hypothesis.  Hence
the left side is positive, so $(9N-M^2)(N-M^2) > 0$.  Since
$M^2 < 3N < 9N$, the first factor is positive.  Hence $N - M^2 > 0$; that is,
$M^2 < N$.
\smallskip

{\em Part (i).}  Taking square roots,
$a/c - b/c = \pm \sqrt{(9N-M^2)(N-M^2)}/(3N-M^2)$.
By the definition of $s$,  $a/c + b/c = s$.
Adding and subtracting gives the formulas in (ii), with one sign for $a/c$
and the other for $b/c$.
\end{proof}

\begin{remark}  Previous versions of this paper gave a proof
of Theorem~\ref{theorem:equilateral3} that used complex numbers
and involved solving an equation using SageMath.  
The proof given above uses only real and hand-checkable algebra.
\end{remark}

\begin{lemma} \label{lemma:ab}
 Let triangle $ABC$ be equilateral.  Suppose it is $N$-tiled using
a tile with angles $(\alpha,\beta,\frac \pi 3)$, where $\alpha$
is not a rational multiple of $\pi$.  Let $(a,b,c)$ be the sides of the 
tile.  Let $M$ be the coloring number of the tiling. Then 
$$ \frac {ab} {c^2} = \frac {4M^2N}{(3N-M^2)^2}. $$
\end{lemma}

\begin{proof}
This is equation (\ref{eq:abc2}) in the proof of Theorem~\ref{theorem:equilateral3}.
\end{proof}

\begin{lemma} \label{lemma:square60}
Let triangle $ABC$ be equilateral.  Suppose it is $N$-tiled using
a tile with angles $(\alpha,\beta,\frac \pi 3)$, where $\alpha$
is not a rational multiple of $\pi$.  Then for some integer
$M < \sqrt{N}$ (the coloring number of the tiling),
$(9N-M^2)(N-M^2)$ is a square.
\end{lemma}

\begin{proof}
Let the sides of the tile be $(a,b,c)$.  According to
Theorem~\ref{theorem:equilateral3}, $M < \sqrt N$, and
the ratios $a/c$ and $b/c$ are given by
\begin{eqnarray*}
  \frac 1 2 \left( \frac {3N+M^2}{3N-M^2}\pm \frac {\sqrt{(9N-M^2)(N-M^2)}}
  {3N-M^2}\right).
\end{eqnarray*}
According to  \cite[Theorem~3.3]{laczkovich2012},
the tile $(a,b,c)$ is rational, so $a/c$ and $b/c$ are rational.
Therefore the expression under the square root is the square of a
rational number; being an integer, it is the square of an integer.
\end{proof}

\begin{theorem} \label{theorem:notprime60}
Let triangle $ABC$ be equilateral.  Suppose it is $N$-tiled using
a tile with angles $(\alpha,\beta,\frac \pi 3)$, where $\alpha$
is not a rational multiple of $\pi$.  Then 
$N$ is not prime.
\end{theorem}

\begin{remark}
(1) In particular
$N$ cannot be 7, 11, or 19.  In \cite{beeson-noseven}, we proved 
that there is no 7 or 11 tiling, and the case studied here was
handled purely computationally, by a method that does not work
for $N=19$.  

(2)  Laczkovich has shown that each rational 
tile with an angle $\gamma = \pi/3$ in which $\alpha$ is not a rational multiple 
of $\pi$ can be used to tile {\em some} equilateral triangle.  But the 
$N$ required might be very large.  Indeed, we tried to construct
a tiling following Laczkovich's instructions in \cite{laczkovich1995},
 but $N$ came out to be over a million, so we could not draw
the tiling.  In 2018, we could not present even one picture of a 
tiling of an equilateral triangle by a tile with an angle $\pi/3$ and incommensurable
angles.  That changed in 2024, as we shall see in the Appendix.  
\end{remark}

\begin{proof}[Proof of Theorem~\ref{theorem:notprime60}] 
Recall the area equation $X^2 = Nab$.
Since the tile is rational  \cite[Theorem~3.3]{laczkovich2012},
 we may re-scale it so that $(a,b,c)$
are integers with no common factor. (That changes the size of $ABC$
and makes $X$ an integer.)  Then $a$, $b$, $c$ are pairwise relatively
prime, since by the law of cosines $c^2 = a^2 - ab + b^2$, a prime dividing
two of them divides the third.  Let $s = a+b+c$.  
Let $M$ be the coloring number of the tiling. 
By Theorem~\ref{theorem:equilateral3}, $M^2 < N$, and by the coloring equation 
$Ms = 3X > 0$, so $M > 0$.
 By the law of cosines,
$$3ab = (a+b)^2 - c^2 = (s-c)^2 - c^2 = s(s-2c).$$
Squaring the coloring equation and using the area equation,
$$M^2s^2 = 9X^2 = 9Nab = 3Ns(s-2c),$$ so $M^2 s = 3N(s-2c)$. 
Multiplying out and isolating the $6Nc$ term on one side, 
\begin{equation}
6Nc = s(3N - M^2). \label{eq:6Nc}
\end{equation}
We claim $s$ and $c$ are relatively prime.  Let $p$ be a prime dividing
$s$ and $c$.  Then $p$ divides $a+b$, and since $c^2 = (a+b)^2 - 3ab$,
$p$ divides $3ab$.  If $p \ne 3$, then $p$ divides $a$ or $b$, and since
$p$ divides $a+b$ it divides both, which is impossible.  If $p = 3$, then
$9$ divides $(a+b)^2 - c^2 = 3ab$, so $3$ divides $a$ or $b$, and again
it divides both.  That proves the claim.  By (\ref{eq:6Nc}), $s$
divides $6N$.  Let $t = 6N/s$, a positive integer.  Then by
(\ref{eq:6Nc}) and the coloring equation,
\begin{eqnarray}
 3N - M^2 &=& tc   \mbox{\  and } \nonumber\\
 t(a+b-c) &=& ts - 2tc \ =\  2M^2. \label{eq:617}
 \end{eqnarray}
Since $s = 6N/t$, we have
$$ X = \frac {M s } 3 = \frac{2MN}{t},$$
and the area equation $X^2 = N ab$ becomes
\begin{equation}
t^2 ab = 4M^2 N. \label{eq:t2ab}
\end{equation}

Now suppose, for proof by contradiction, that $N$ is prime.  If $N$ is
2 or 3, then $M = 1$, since $0 < M < \sqrt N$, and $(9N-M^2)(N-M^2)$ is
17 or 52, not a square, contradicting Lemma~\ref{lemma:square60}.  So $N \ge 5$.  Since
$a+b-c$ is a positive integer, $t \le t(a+b-c) = 2M^2 < 2N$.  If $t = N$,
then by \eqref{eq:617}, $N(a+b-c) = 2M^2 < 2N$, so $a + b - c = 1$ and $N = 2M^2$, which is
even.  The divisors of $6N$ less than $2N$, other than $N$, are
1, 2, 3, and 6.  Since $t = 6N/s$ is an integer,  $t$ is one of these, and $N$ does not divide $t$.
By (\ref{eq:t2ab}), $N$ divides $ab$.  By symmetry we may suppose
$N$ divides $a$, say $a = Na'$.  By the triangle inequality, $2a < s = 6N/t$,
so $a't < 3$.  There are three cases.

Case 1, $t = 1$ and $a' = 1$.  Then $a = N$, and (\ref{eq:t2ab}) gives
$b = 4M^2$.  By \eqref{eq:617} with $t = 1$, $a + b - c = 2M^2$, so
\begin{eqnarray*}
c =  N + 2M^2  && \mbox{\qquad since $a=N$ and $b = 4M^2$} \\
6N = s &&\mbox{\qquad since $t = 6N/s = 1$}\\
 6N  = a + b + c &&\mbox{\qquad by definition of $s$}\\
 6N  = N + 4M^2 + N + 2M^2 &&\mbox{\qquad substituting $a$, $b$, $c$}\\
6N   = 2N + 6M^2  &&\\
2N = 3M^2  &&\mbox{\qquad subtracting $2N$ and dividing by 2}
\end{eqnarray*}
Then $3$ divides $N$, so $N = 3$, but $N \ge 5$.
\smallskip

Case 2, $t = 2$ and $a' = 1$.  Then $a = N$, and
\begin{eqnarray*}
4Nb  = 4M^2N              &&\mbox{\qquad by (\ref{eq:t2ab}), since $t = 2$ and $a = N$}\\
b  =  M^2                &&\mbox{\qquad dividing by $4N$}\\
2(a+b-c)  =  2M^2         &&\mbox{\qquad by (\ref{eq:617}), since $t = 2$}\\
c  = a + b - M^2 \     &&\mbox{\qquad dividing by 2 and solving for $c$}\\
c = a                  &&\mbox{\qquad since $b = M^2$}\\
a^2  =  a^2 - ab + b^2  &&\mbox{\qquad law of cosines, since $c = a$}\\
b  =  a               &&\mbox{\qquad since $b(b - a) = 0$ and $b > 0$}
\end{eqnarray*}
So the tile is equilateral and $\alpha = \pi/3$, contradicting the hypothesis.
\smallskip

Case 3, $t = 1$ and $a' = 2$.  Then $a = 2N$, and
\begin{eqnarray*}
2Nb = 4M^2N  &&\mbox{\qquad by (\ref{eq:t2ab}), since $t = 1$ and $a = 2N$}\\
b = 2M^2   &&\mbox{\qquad dividing by $2N$}
\end{eqnarray*}
So $a$ and $b$ are both even, contradicting their relative primality.

All cases lead to a contradiction, so $N$ is not prime.

\end{proof}

\subsection{An algorithm to decide if there is an $N$-tiling with $\gamma = \pi/3$}
 
\begin{theorem}  \label{theorem:decidable}
Given $N$,  
there is a finite set $S$ of  not more than $\sqrt {N}$ tiles 
such that, if any tile with angles $(\alpha,\beta,\pi/3)$ and 
$\alpha$ not a rational multiple of $\pi$ can $N$-tile an equilateral
triangle, then 
one of the tiles in $S$ can do so.  
\end{theorem}

\begin{proof} By Theorem~\ref{theorem:equilateral3},
the coloring number $M$ of any $N$-tiling of an equilateral triangle 
$ABC$ is at most $\sqrt{N}$,
and $(N,M)$ together determine $a/c$ and $b/c$, 
which are rational by Lemma~\ref{lemma:square60}; hence they 
determine the sides $(a,b,c)$ as integers with no common factor. 
This provides the finite set $S$ of tiles, and Theorem~\ref{theorem:equilateral3}
says that any $N$-tiling uses one of those tiles.  The triangle $ABC$
will be scaled by the area equation so that it has the same area as $N$ tiles.
\end{proof}

\begin{remark}  It follows that whether such a tiling
exists is computable in a finite (though perhaps large) number 
of steps, using standard search algorithms. We have not 
reported on such computer searches in this paper, but
we briefly discuss the problem.  Given $N$, we determine
the finite set $S$ of possible tiles and triangle $ABC$,
and test one by one whether
that tile can tile that $ABC$.  Thus the problem reduces to 
testing whether a specific tile can tile a specific triangle.

We can consider connected partial tilings as nodes in a graph, where
there is an edge between two nodes $p$ and $q$ if partial tiling $q$
is obtained from $p$ by adding one more tile within the boundaries
of $ABC$, the new tile sharing at least one vertex and at least part
of at least one edge with the tiling $p$, and not overlapping any
tile of $p$.   Each partial tiling $p$
has finitely many neighbors, which can be algorithmically generated
by enumerating the possible ways to extend a given partial tiling.
In so doing we need to test whether a triangle overlaps another triangle;
the precision issue involved is discussed below.  
  No path has length more than $N$, 
since the area of a partial tiling cannot exceed $N$ times the area
of the tile.  Finally, we can test algorithmically whether a 
given partial tiling is actually a tiling of $ABC$; again a precision 
issue arises.  The point of these precision issues is that 
we need to determine
in a finite number of steps of computation whether two points 
coincide or not, whether two tile edges coincide or not, and whether 
a point lies on a given line segment.  Standard graph-search 
algorithms can then be used to search for a tiling.
   
The use of such programs is an active area of research,  but we 
do not report on it here, because of the twin difficulties of 
writing such programs correctly, and certifying their correctness,
so that if the program finds,  for example, no $60$-tiling 
of an equilateral $ABC$ by $(3,5,7)$ is possible, we could obtain
a formal proof of that result.   
\end{remark}

\subsection{Comparison with Laczkovich's results}

Given $N$, we have shown how to determine a finite set of 
``possible tiles'' that include every tile (with a $\pi/3$ angle
and $\alpha$ not a rational multiple of $\pi$)  that can be used
in an $N$-tiling of an equilateral triangle.  In other words,
given $N$, the tile (and hence $ABC$) are determined; or more 
accurately, all other possibilities for $(N,ABC)$ are eliminated.
We do not know if a tiling really exists, except by trial and error.
 
 Laczkovich 
proved that any such tile can be used to $N$-tile an equilateral
triangle $ABC$, if we choose $N$ large enough.  In other words,
given the tile,  at least one pair $(N, ABC)$ is determined such
that a tiling exists.   

Comparing those results, the question naturally arises, whether
the tile actually determines $N$.   That is, can the same 
tile be used for tiling two  equilateral triangles of different sizes?
Well, given one tiling, we can always replace each tile by a 
quadratic tiling of $m^2$ smaller tiles, thus producing an $m^2N$ 
tiling.  So the best we can hope for is that the squarefree part
of $N$ might be determined by the tile.
Since we have explicit formulas for $a/c$ and $b/c$, this question
can be answered.

\begin{theorem} \label{theorem:comparison}
Suppose $T$ is a triangle with angles $(\alpha,\beta,\pi/3)$
with $\beta$ not a rational multiple of $\pi$.  
Suppose $T$ can be used to $N$-tile an 
equilateral triangle, and also to $K$-tile a (different) equilateral triangle.
Then $N$ and $K$ have the same squarefree part.  
\end{theorem}

\begin{proof}
Let $M$ be the coloring number of the $N$-tiling and $J$ the
coloring number of the $K$-tiling.
By (\ref{eq:s}) applied to both tilings, 
$$\frac{a+b} c = \frac{3N+M^2}{3N-M^2} $$
and 
$$ \frac{a+b} c= \frac{3K+J^2}{3K-J^2}.$$
Therefore
\begin{eqnarray}
\frac {3N+M^2}{3N-M^2}  
&=& \frac {3K+J^2}{3K-J^2}   \label{eq:689}
\end{eqnarray}
Cross-multiplying, we have
$$ (3K-J^2)(3N+M^2) = (3N-M^2)(3K+J^2).$$
Simplifying, 
$$ J^2N = M^2 K.$$
  Hence $N$ and $K$ have the same square-free part.
\end{proof}

\section{Laczkovich's graph $\Gamma_c$} 

Laczkovich has proved  (already in 1995 \cite{laczkovich1995}) that, given a rational 
tile of the shape we are now considering, there is an $N$-tiling of 
{\em some} sufficiently large equilateral $ABC$; but $N$ might have
to be large.  As mentioned above, following Laczkovich we 
found $N = 10935$,  which Herdt improved in 2024 to $N=1215$.
 See Figs.~\ref{figure:10935} and ~\ref{figure:1215}. 

In 2012 \cite{laczkovich2012}, Laczkovich 
made a significant advance:  a tile with a $2\pi/3$ angle 
that tiles an equilateral triangle must either have both 
its other angles $\pi/6$, or else both the following conditions 
hold: 
\smallskip

 (i) the tile is rational (that is, the ratios
of the sides are rational), and 
\smallskip

(ii) the other two angles of the tile are not rational multiples of $\pi$.
\smallskip

These statements are proved in  Theorem~3.3  and Lemma~3.2 of \cite{laczkovich2012}, respectively.  It is the rationality of 
the tile that is the significant advance of 2012 \cite{laczkovich2012}, as (ii) was already 
proved in 1995 \cite{laczkovich1995}.
The main tool is a directed graph $\Gamma_c$.  We give the definition and slight 
modifications of two important lemmas, which we will apply below.

We need some terminology.  Given a tiling of (in our case) a triangle $ABC$,
an {\em internal segment} is a line segment connecting two vertices
of the tiling that is contained in the union of the boundaries of the 
tiles, and lies in the interior of $ABC$ except possibly for its endpoints.
A {\em maximal segment} is an internal segment that is not part of a longer
internal segment.  A {\em left-maximal segment} is an internal segment 
$XY$ that is not contained in a longer segment $UXY$, i.e., a segment 
$UY$ with $X$ between $U$ and $Y$.  A tile is {\em supported by}
$XY$ if one edge of the tile lies on $XY$.   The internal segment $XY$ is 
said to have ``all $c$'s on the left'' if the endpoints $X$ and $Y$ are
vertices of tiles supported by $XY$ and lying on the left side of $XY$,
and all tiles supported by $XY$ lying on the left of $XY$ have there 
$c$ edges on $XY$.  Similarly for ``all $c$'s on the right.'' 

An internal segment $XY$ is said to {\em witness the relation} 
$jc = \ell a + mb$ in case the endpoints
$X$ and $Y$ are vertices of tiles on both sides of $XY$, and either  
\begin{itemize}
\item  $XY$ has all $c$'s on one side, and exactly $j$ of them (that 
is, the length of $XY$ is $jc$), and on the other side $XY$ supports
$\ell$ tiles with their $a$ edges on $XY$ and $m$ tiles with their
$b$ edges on $XY$ (in any order) and no other tiles, or
\item  $XY$ has $j+n$ tiles with $c$ edges on one side, and on 
the other side $X$ supports $\ell$ tiles with their $a$ edges on $XY$ and $m$ tiles with their
$b$ edges on $XY$ and $n$ tiles with their $c$ edges on $XY$.
\end{itemize}

\noindent{\em Exercise}.  Identify the relations that are witnessed in 
Fig.~\ref{figure:bryce1440}.
\medskip

This condition implies 
that $c$ is not a linear combination of $a$ and $b$ with nonnegative
rational coefficients, but it is stronger than that statement, in some
way limiting the size of the (numerators and denominators of the) coefficients.
In particular, if $XY$ witnesses $jc = \ell a + m b$, then since 
$XY$ is an interior segment, the length of $XY$ is less than 
the diameter (which for an equilateral triangle is the side)  of $ABC$.  That
places a bound on $\ell$ and $m$.

Similarly we use the terminology ``$XY$ witnesses a relation 
$jc = \ell a + mc$'',  which implies ``$c$ is a rational multiple of $a$'',
but is stronger.  

Laczkovich defined a tiling to be {\em regular} if there are two angles
(say $\alpha$ and $\beta$) of the tile such that at each  vertex $V$
of the tiling, the number of tiles having angle $\alpha$ at $V$ is the 
same as the number of tiles having angle $\beta$ at $V$.  According to 
Lemma~3.2 of \cite{laczkovich2012}, in an irregular tiling, 
$(\alpha,\beta,\gamma)$ are linear combinations with rational coefficients
of the angles of the tiled polygon.  In this paper the tiled polygon is 
the equilateral triangle, so an irregular tiling has angles that are 
rational multiples of $\pi$.  Therefore, the case of interest in this 
section, when $\alpha$ and $\beta$ are not rational multiples of $\pi$ and 
$\gamma = 2\pi/3$, only can occur in regular tilings.  

In the work below, we shall make use of Lemmas~4.5 and 4.6 of 
\cite{laczkovich2012}.  These lemmas make use of the directed
graph $\Gamma_c$ defined on page 346.  We review that definition 
now.

\begin{definition} [The directed graph $\Gamma_c$] Given a tiling of some triangle,
the nodes of $\Gamma_c$ are certain vertices of the tiling. A {\em link} of
$\Gamma_c$  connects $X$ to $Y$ if there is a maximal segment $XZ$ containing $Y$, with
tiles on both sides of it (so no tile edge on either side extends $XZ$ past $X$),
and the two sides can be called Side~1 and Side~2 so that
(i) $XY$ is covered by $c$ edges on Side~1;
(ii) past $Y$, on Side~1, there is a tile with a vertex at $Y$ whose edge on $XZ$
     is not a $c$ edge;
(iii) the tile on Side~2 with a vertex at $X$ does not have its $c$ edge on $XZ$.
$\Gamma_a$ and $\Gamma_b$ are defined similarly.
This is a directed graph: $X$ is the {\em tail} of the link and $Y$ is the {\em head}.
\end{definition}

\begin{remark} This is Laczkovich's definition \cite{laczkovich2012}, p.~346.
The graph $\Gamma_a$ of \cite{beeson-triangletiling3} is slightly different: it omits (iii), and
instead assumes that $Y$ lies in the interior of an edge on Side~2, which
Laczkovich derives.   We use the word {\em link} for the edges of $\Gamma_c$ since {\em edge} is 
already in use for tile edges.
\end{remark}

\noindent
{\em Example}. In Fig.~\ref{figure:10935}, look at the longest side 
of one of the light blue components of the tiling.  That segment is 
composed of 21 $c$ edges.  At one end it cannot be extended:  that is $X$.
At the other end, it does extend beyond the blue triangles, but there it 
has $a$ edges on both sides.  At $X$, the tile on the other side
from the $c$ edges has its $b$ edge on $XY$, meeting condition~(iii).
Hence, there is a link of $\Gamma_c$ from 
$X$ to $Y$.  
\medskip

\noindent
{\em Exercise}. Identify the graphs $\Gamma_a$, $\Gamma_b$, and $\Gamma_c$ in 
Fig.~\ref{figure:bryce1440} and  Fig.~\ref{figure:10935}.
\medskip

We now state the versions of Laczkovich's lemmas that we need.
These lemmas presuppose a regular 
tiling of a convex polygon, in our application 
an equilateral triangle.  

\begin{lemma}[Laczkovich's Lemma~4.5] \label{lemma:propagation2}
Suppose the tiling does
not witness any relation $jc = \ell a + mb$.  Let $XY$ be 
a segment of the tiling (internal or on the boundary) and let 
$V$ be a vertex of the tiling lying on the interior of 
$XY$, lying either on the boundary
or on an internal point of an edge of some tile. Suppose that of the 
two tiles supported by $XY$ with a vertex at $V$,  one has edge
$c$ on $XY$ and the other has edge $a$ or $b$ on $XY$.  Then 
there is a link of $\Gamma_c$ starting from $V$. 
\end{lemma} 

\noindent{\em Remark}. In Fig.~\ref{figure:10935},  the edge
of $\Gamma_c$ mentioned as an example does witness a relation
$21 c = 24 a + 15b$.    
\medskip

\begin{proof} Laczkovich's statement replaces the first 
sentence of the lemma by ``$c$ is not a linear combination of $a$ 
and $b$ with nonnegative rational coefficients.''  But the proof 
actually proves our version; only the last sentence uses the 
weakened hypothesis. \footnote{
Laczkovich's proof is easier in
case $ABC$ is equilateral, which is the case we need.  Namely,
 if $\gamma = \pi/3$
then the complicated second paragraph is not needed, and if $\gamma = 2\pi/3$
then in that paragraph we have $p=3$.}
\end{proof}

\begin{lemma}\label{lemma:noboundarylink}
Suppose there is a tiling of an equilateral triangle by a tile with 
angles $(\alpha,\beta,\gamma)$ and sides $(a,b,c)$, where
 $\gamma = \pi/3$ or $2\pi/3$ and $\alpha$ not a rational multiple of $\pi$.
If the tiling does not witness any relation 
of the form $jc = \ell a + mb$,  then the in-degree and out-degree are
equal at every vertex of $\Gamma_c$, and 
 no link of $\Gamma_c$ starts on the boundary of $ABC$.
\end{lemma}

\begin{proof}
If the tiling does not witness any relation $jc = \ell a + mb$,
then by Lemma~\ref{lemma:propagation2}, each ending point (head) of an edge is the starting point (tail)
of a new edge.  At most one link ends at any vertex, so the out-degree of each node is at least the in-degree.
Since every edge has a head and a tail, the total in-degree is equal to 
the total out-degree.  Therefore the out-degree at each node is equal to the in-degree.
In particular, by 
the definition of $\Gamma_c$, no ending point of a link lies on the boundary
of $ABC$; hence no link begins on the boundary.  
\end{proof}

\section{A tile $(\alpha,\beta,2\pi/3)$ with $\alpha/\pi$ irrational}

Every vertex of the tiling with total angle $\pi$ either is 
composed of $\alpha+\beta+\gamma$ or of $3\alpha+3\beta$.  
Since the latter form has six tiles meeting at the vertex, there
is no coloring equation, since that would require an 
odd number of tiles at each vertex with total angle $\pi$. Even in a tiling without 
such vertices, there still could not be a coloring equation, 
because there will have to be a ``center'' somewhere in the tiling, with three tiles each
having angle $2\pi/3$.  The existence of a ``center'' follows
from the observation that since $\alpha+\beta = \pi/3$, at each 
vertex of $ABC$, there will have to be two tiles with angles 
$\alpha$ and $\beta$; we do not give details since our only purpose
here is to explain why we cannot use the coloring equation for these
tilings. 

Before coming to the main theorem, we prove a lemma.
It may seem obvious, but it does actually need 
a proof.
\begin{lemma} \label{lemma:notallsupported}
 In an $N$-tiling of any triangle $ABC$ with $N > 3$, no segment of the tiling
can support all $N$ tiles.
\end{lemma}

\noindent{\em Proof}.  Suppose segment $XY$ supports all $N$ tiles.
If there is no vertex of the tiling lying on the interior of $XY$ then 
there is only one tile on each side of $XY$.  Then $N=1$ if $XY$ 
lies on the boundary of $ABC$, and otherwise $N=2$, contradiction.
Therefore there is a vertex $V$ between $X$ and $Y$. If at least 
three tiles meet at $V$ on the same side of $XY$, then the middle 
one is not supported by $XY$.  Hence exactly two tiles on the same
side of $XY$ meet at $V$.  Then both tiles must have their $\gamma$
angle at $V$,
 since otherwise the sum of the angles cannot be $\pi$.
Then $\gamma = \pi/2$.  Let $T$ be one of those two tiles and 
let $W$ be its other vertex on $XY$.  Then $T$ does not have its 
$\gamma$ angle at $W$.  If $W$ is in the interior of $XY$ 
then the other two tiles on the same side of $XY$ as $T$ and 
with a vertex at $W$ must together make more 
than a $\pi/2$ angle at $W$, so there must be at least two such tiles.
Then the one adjoining $T$ is not supported by $XY$, contradiction. 
Therefore $W$ is not in the interior of $XY$, but must be $X$ or $Y$.
Hence only two tiles on that side of $XY$ are supported by $XY$,
and their other sides (the ones not shared between the two tiles 
or lying on $XY$) lie on the boundary of $ABC$, and the angles 
of those tiles at $X$ and $Y$ are acute.  Now if $XY$ lies on 
the boundary of $ABC$, then $N=2$, contradiction.  Hence $XY$ is an 
interior segment.  Then the same argument applies to the other side 
of $XY$, so there are exactly two tiles supported on that side of 
$XY$ as well, whose other edges lie on the boundary of $ABC$ and 
have acute angles at $X$ and $Y$. But then $ABC$ is a
 quadrilateral, with diagonal $XY$.  (Note that some quadrilaterals
can indeed be tiled in such a way that the diagonal supports all four 
tiles.)   That contradicts the hypothesis that $ABC$ is a triangle.
That completes the proof.

\begin{theorem}
\label{theorem:notprime120}
Let equilateral triangle $ABC$ be $N$-tiled by a 
tile with angles $(\alpha,\beta,\gamma)$, with 
$\gamma = 2\pi/3$ and $\alpha$ not a rational multiple of $\pi$.
Then $N$ is not a prime number.
\end{theorem}

\begin{proof}
By Theorem~3.3 of \cite{laczkovich2012}, $(a,b,c)$ are pairwise
commensurable.  Without loss of generality we can
assume that $(a,b,c)$ are integers with no common factor. 
As explained above, Lemma~3.2 of \cite{laczkovich2012} implies 
the tiling is regular.  
Assume, for proof by contradiction, that $N$ is a prime number.

Let $X$ be the length of the sides of equilateral $ABC$.  
Then the area equation is 
\begin{eqnarray*}
X^2 \sin(\pi/3) &=& Nab \sin(2\pi/3)
\end{eqnarray*}
Since $\sin(\pi/3) = \sin(2\pi/3)$ we have
\begin{eqnarray*}
X^2 &=& Nab  \mbox{\qquad area equation} 
\end{eqnarray*}
If $N \le 3$, then $X^2 = Nab  \le 3ab < a^2 + ab + b^2 = c^2$.
So $c > X$. But a tile’s $c$ edge lies inside $ABC$, so $c \le X$, 
since the diameter of an equilateral triangle is its side. 
Now we have $c > X$ and $c \le X$, a contradiction.

Since each side of $ABC$ is the disjoint union of a set 
of tile edges, we have for some non-negative integers $(p,q,r)$,
$$ X = pa + qb + rc.$$
Then $X$ is an integer.   Since $X^2 = Nab$,  $N$ divides $X^2$,
which we write as usual $N \vert X^2$.  Then we have
\begin{eqnarray*}
&& N \mid X^2 \\
&& N \mid X  \mbox{\qquad since $N$ is prime} \\
&& N^2 \mid X^2 \\
&& N^2 \mid Nab  \mbox{\qquad since $X^2 = Nab$} \\
&& N \mid ab
\end{eqnarray*}
Since $N$ is prime,  $N$ divides $a$ or $N$ divides $b$.
Since so far, nothing distinguishes $\alpha$ from $\beta$ except 
the name, we may assume without loss of generality that $N$ divides $a$.
Then there is an integer $e \ge 0$ such that $a = Ne$.
 Since $q$ and $r$ are counts of tiles 
supported by $AB$,  $q + r \le N$; and
by Lemma~\ref{lemma:notallsupported}, $q+r < N$.  
 
By the law of cosines, we have 
\begin{eqnarray*}
c^2 &=& a^2 + b^2 - 2ab \cos(2\pi/3) \\
c^2 &=& a^2 + b^2 + ab  \mbox{\qquad since $\cos(2\pi/3)  = -1/2$}  \\
c^2 &=& N^2 e^2 + b^2 + Neb 
\end{eqnarray*}
Therefore  $c$ is congruent to $\pm b \pmod N$ and 
also $N$ does not divide $b$, since if $N \mid b$ then also $N \mid c^2$ 
and hence $N \mid c$, contradiction, since then $N$ would divide
all three of $(a,b,c)$, but $(a,b,c)$ have no 
common divisor.   
We will show that $AB$ supports only tiles
with their $a$ or $b$ edges on $AB$; that is, no $c$ edges.  We argue by cases.
\smallskip

{\em Case~1}, $c \equiv b \pmod N$.  Then $X \equiv pa + (q+r)b \pmod N$.
Since $N \mid a$, $X \equiv (q+r) b \pmod N$.   Since $N \nmid b$,
and $N$ is prime, we have $N \mid (q+r)$.  Since $q+r < N$, 
we have $q+r = 0$. Since $q$ and $r$ are nonnegative integers, we have 
$q = r = 0$.   Therefore,  every tile supported by $AB$ has
its $a$ edge on $AB$.   In particular, no $c$ edges, as claimed.
\smallskip

{\em Case~2}, $c \equiv -b \pmod N$. 
We will show that 
no relation $jc = \ell a + m b$  with $j > 0$ is witnessed.  Suppose it is; 
since $c \equiv - b \pmod N$, we have
then  $- jb \equiv \ell a + mb \pmod N$, so $N \mid m + j$.  But $m$ and $j$
are counts of tiles on opposite sides of the witnessing segment,
so $m+j \le N$.   
Since $N \mid m+j$,  and $j > 0$, we must have $m+j = N$, and every tile is supported
by the same segment, contradicting Lemma~\ref{lemma:notallsupported}.
Hence, as claimed, no relation $jc = \ell a + m b$ is witnessed. 
Then by Lemma~\ref{lemma:noboundarylink},
no link of $\Gamma_c$ starts on the boundary.   However, by Lemma~\ref{lemma:propagation2},
a link starts on the boundary at any vertex between a $c$-edge and a non-$c$ edge.
Therefore, there is no vertex on the boundary of $ABC$ with a $c$ edge on one side
and a non-$c$ edge on the other side.   Therefore, each side is either all $c$-edges or 
no $c$-edges.

Since $X = pa + qb + rc$, we have
$$X \equiv (q-r)b \pmod N.$$
Therefore $N \mid (q-r)$.  We have $\vert q-r \vert  \le q+r <  N$.
 Hence $q-r = 0$; that is, $q=r$. 
Hence $AB$ is not all $c$-edges, since that would mean $q = 0 < r$. 
Since each side is all $c$-edges or no $c$-edges, $AB$ has no $c$ edges, as claimed.
That concludes the argument by cases to demonstrate that point.

Then each tile supported by 
$AB$ has a $\gamma$ angle on $AB$ (since the $c$ edge of each tile is in the interior).
 Since $\gamma = 2\pi/3 > \pi/2$, there cannot be
two  $\gamma$ angles at the same vertex,  and there cannot be a $\gamma$ angle
at the vertices $A$ and $B$.  Then some number $k$ of $\gamma$ angles 
are distributed over $k-1$ vertices, at most one at each vertex.  
That contradicts the pigeonhole principle.
\end{proof}

\section{$ab$ divides $X$ and $N$}
The following theorem and its proof are due to J. P. Harries \cite{harries2026, harries2026no19tiling}, who kindly allowed us to include it.
It confirms Zhang’s Conjecture~1 \cite{zhang2026}, and the corresponding part of his Conjecture~3 (the case $\gamma = \pi/3$) \cite{zhang2026}.
It also gives a new proof that $N$ is not prime.

\begin{theorem} \label{theorem:abdividesX}  Suppose the equilateral triangle $ABC$ is tiled by 
a tile with angles $(\alpha,\beta,\gamma)$, where $\alpha$ is not a rational multiple of $\pi$,
and $\gamma = \pi/3$ or $2\pi/3$.
Suppose the tile has integer side lengths $(a,b,c)$ with no common factor.   Then 
$ab$ divides the side length $X$ of $ABC$. 
\end{theorem}

\begin{proof}  The argument is due to Harries {\em op.~cit.} It uses
Laczkovich's method of an additive function on directed segments
\cite{laczkovich2012}.  Since $\gamma \in \{\pi/3, 2\pi/3\}$ and $\alpha$ is
not a rational multiple of $\pi$, neither is $\beta = \pi - \gamma - \alpha$.

{\em Directions.}  Place $AB$ on the $x$-axis, with $ABC$ above it, and let
$G = \{ j\pi/3 + k\alpha : j, k \in \Z\}$, an additive group of angles
taken mod $2\pi$.  The angles of the tile are $\alpha$, $\gamma$, and
$\beta = \pi - \gamma - \alpha$, all in $G$.  Hence if one edge of a tile has
direction in $G$, so do its other two edges.  The tiles supported by $AB$
have an edge in direction $0$, and any two tiles sharing a segment of
positive length have parallel edges, whose directions differ by $0$ or $\pi$.
Since the union of the tiles is connected, every edge of every tile has its
direction in $G$.  The representation $\theta = j\pi/3 + k\alpha$ is unique,
with $j$ determined mod $6$: if $j\pi/3 + k\alpha \equiv j'\pi/3 + k'\alpha
\pmod{2\pi}$, then $(k-k')\alpha$ is a rational multiple of $\pi$, so $k = k'$,
and then $j \equiv j' \pmod 6$.

{\em The function $\Phi$.}  Define $\chi : G \to \Z[z,z^{-1}]$ by
$\chi(j\pi/3 + k\alpha) = (-1)^j z^k$.  This is well defined, and
$\chi(\theta + \varphi) = \chi(\theta)\chi(\varphi)$; in particular
$\chi(\theta + \pi) = -\chi(\theta)$.  For a directed segment $PQ$ with
direction $\theta \in G$, put $\Phi(PQ) = |PQ|\,\chi(\theta) \in
\R[z,z^{-1}]$, and for a polygon $R$ whose sides have directions in $G$, let
$\Phi(\partial R)$ be the sum of $\Phi$ over the sides of $R$, traversed
counterclockwise.  Then
\begin{equation}
\Phi(\partial ABC) = \sum_i \Phi(\partial T_i), \label{eq:Phisum}
\end{equation}
where $T_i$ are the tiles.  Indeed, each maximal segment in the interior of
$ABC$ is covered on its two sides by tile edges of the same total length,
traversed in opposite directions, so their contributions cancel since
$\chi(\theta+\pi) = -\chi(\theta)$; and the tile edges on the boundary of
$ABC$ add up to $\Phi(\partial ABC)$.

{\em One tile.}  Let $\sigma = 1$ and $d = a - b$ if $\gamma = 2\pi/3$, and
$\sigma = -1$ and $d = a + b$ if $\gamma = \pi/3$.  By the law of cosines,
$c^2 = a^2 + \sigma ab + b^2$, so in both cases
\begin{equation}
c^2 - d^2 = 3\sigma ab. \label{eq:cdab}
\end{equation}
Consider a tile with its $a$ edge traversed in direction $0$, the tile lying
above it.  The ends of the $a$ edge carry the angles $\beta$ and $\gamma$.  If
$\gamma$ is at the east end, the next two sides are $b$, in direction
$\pi - \gamma$, and $c$, in direction $2\pi - \gamma - \alpha$; this gives
$\Phi = a + b\,\chi(\pi-\gamma) + c\,\chi(2\pi - \gamma - \alpha)
= d + \sigma c z^{-1}$.  If $\beta$ is at the east end, the next two sides are
$c$, in direction $\pi - \beta = \gamma + \alpha$, and $b$, in direction
$2\pi - \beta - \alpha = \pi + \gamma$; this gives
$\Phi = a + c\,\chi(\gamma + \alpha) + b\,\chi(\pi + \gamma)
= d + \sigma c z$.  (Here $\chi(\pi/3) = -1$, $\chi(2\pi/3) = 1$.)  Every tile
of the tiling, reflected or not, is obtained from one of these two by a
translation, which does not change $\Phi$, and a rotation by some
$\theta \in G$, which multiplies $\Phi$ by $\chi(\theta) = \pm z^k$.
Therefore each $\Phi(\partial T_i)$ has the form
$\pm z^{k}(d + \sigma c z^{\pm 1})$.

{\em The boundary.}  The sides of $ABC$, traversed counterclockwise, have
directions $0$, $2\pi/3$, $4\pi/3$, and $\chi$ is $1$ at each of them.  So
$\Phi(\partial ABC) = 3X$, and by (\ref{eq:Phisum})
\begin{equation}
3X = \sum_i \varepsilon_i z^{k_i} (d + \sigma c z^{\eta_i}),
\qquad \varepsilon_i, \eta_i \in \{1,-1\},\ k_i \in \Z. \label{eq:3X}
\end{equation}
The right side has integer coefficients, so (\ref{eq:3X}) is an identity in
$\Z[z,z^{-1}]$.

{\em Reduction mod $3ab$.}  First, $a$, $b$, $c$ are pairwise relatively
prime: a prime dividing two of them divides the square of the third, by
$c^2 = a^2 + \sigma ab + b^2$, contradicting the absence of a common factor.
Next, $3 \nmid c$: otherwise $d^2 = c^2 - 3\sigma ab \equiv 0 \pmod 3$, so
$3 \mid d$, so $9 \mid c^2 - d^2 = 3\sigma ab$ and $3$ divides $a$ or $b$;
then $3 \mid d = a \mp b$ gives $3 \mid a$ and $3 \mid b$, which is
impossible.  Hence $c$ is invertible mod $3ab$.  Also $d$ is invertible
mod $3ab$: a prime $p \ne 3$ dividing $ab$ divides exactly one of $a$, $b$,
hence not $d = a \mp b$; and $3 \nmid d$, since $d^2 \equiv c^2
\not\equiv 0 \pmod 3$.  Let $u = -\sigma d c^{-1}$, a unit mod $3ab$, and let
$f : \Z[z,z^{-1}] \to \Z/3ab\Z$ be the ring homomorphism with $f(z) = u$.
Then $f(d + \sigma c z) = d - d = 0$ and, using (\ref{eq:cdab}),
$$f(d + \sigma c z^{-1}) = d - c^2 d^{-1} = d^{-1}(d^2 - c^2)
= -3\sigma ab\, d^{-1} = 0.$$
Applying $f$ to (\ref{eq:3X}) gives $3X \equiv 0 \pmod{3ab}$.  Hence
$ab \mid X$.
\end{proof}

Theorem~\ref{theorem:abdividesX} gives a new proof that $N$ is not 
prime.

\begin{corollary} \label{corollary:abdividesN}
Suppose the equilateral triangle $ABC$ is tiled by 
a tile with angles $(\alpha,\beta,\gamma)$, where $\alpha$ is not a rational multiple of $\pi$,
and $\gamma = \pi/3$ or $2\pi/3$.
Suppose the tile has integer side lengths $(a,b,c)$ with no common factor.   Then 
$N = m^2 ab$  for some $m$.   In particular, $N$ is not prime.
\end{corollary}

\begin{proof} The area equation is $X^2 \sin(\pi/3) = Nab \sin \gamma$, where $\sin \gamma = \onehalf \sqrt 3$ is the same
for $\gamma = 2\pi/3$ and $\pi/3$.  Then $X^2 = Nab$.  By Theorem~\ref{theorem:abdividesX}, there is a positive integer $m$ 
such that $X = mab$.
Then $m^2(ab)^2 = N ab$,  so $N = m^2 ab$.   We claim that neither $a$ nor $b$ can be 1.
 If $a = 1$ and $b \ge 2$, then $c^2 = b^2 \pm b + 1$ lies strictly between $b^2$ and $(b+1)^2$, 
 or between $(b-1)^2$ and $b^2$, so it is not a square; similarly if $b = 1$. And $a = b = 1$ would give $\alpha = \beta$, 
 impossible since $\alpha$ is not a rational multiple of $\pi$. Therefore neither $a$ nor $b$ is 1, as claimed. 
 Therefore $N$ is not prime.
\end{proof}

\section{Necessary conditions}
%
%
%
 In this 
section we restrict attention to tiles with angles $(\alpha,\beta,\pi/3)$
where $\alpha$ and $\beta$ are not  rational multiples of $\pi$.  Suppose
there is an $N$-tiling of an equilateral triangle: what can we say 
about $N$?  Above we proved that $N$ cannot be prime; but certainly
there are many composite $N$ for which there is no $N$-tiling.  
Let us consider the question, given $N$, can we determine whether or 
not there is an $N$-tiling of an equilateral triangle with $\gamma = \pi/3$
and $\alpha$ not a rational multiple of $\pi$?  We do have one important
necessary condition already:  $ab$ must divide $N$, by Corollary~\ref{corollary:abdividesN}.
 
We have reduced the problem to three computational steps:

\begin{itemize}
\item  Determine if the equations of Theorem~\ref{theorem:equilateral3}, part~(i)
 have
rational solutions for the ratios $a/c$ and $b/c$.  If they do 
not, there is no $N$-tiling.  If they do, let $(a,b,c)$ be 
integers with no common divisor whose ratios solve the equations.

\item  Reject the candidate $(a,b,c)$ if $ab$ does not divide $N$.

\item  Determine, given $N$ and $(a,b,c)$   whether
a tiling actually exists.
\end{itemize}

We shall not undertake the third step,  but concentrate on the first two.
The first step can be done for $\gamma = \pi/3$ by using the coloring number $M$,
since $M$ is bounded in 
terms of $N$, and $M$ and $N$ together determine the tile.
The results are shown in Table~2.

\begin{table} [ht]
\caption{Tilings with $\gamma = \pi/3$  not ruled out by the area and \\
coloring equations and Theorem~\ref{theorem:abdividesX}.}
\label{table:2}
\begin{center}
\begin{tabular}{rrrr}
$N$ &  $M$ & $(a,b,c)$ & side of $ABC$ \\
\hline
24 & 4 & (3, 8, 7) & 24\rlap{\enspace\nomark} \\
40 & 6 & (5, 8, 7) & 40 \\
96 & 8 & (3, 8, 7) & 48 \\
105 & 7 & (5, 21, 19) & 105 \\
105 & 9 & (7, 15, 13) & 105 \\
120 & 10 & (8, 15, 13) & 120 \\
160 & 12 & (5, 8, 7) & 80 \\
216 & 12 & (3, 8, 7) & 72 \\
280 & 10 & (7, 40, 37) & 280 \\
336 & 18 & (16, 21, 19) & 336 \\
360 & 18 & (5, 8, 7) & 120 \\
384 & 16 & (3, 8, 7) & 96 \\
385 & 15 & (11, 35, 31) & 385 \\
420 & 14 & (5, 21, 19) & 210 \\
420 & 18 & (7, 15, 13) & 210 \\
480 & 20 & (8, 15, 13) & 240 \\
585 & 13 & (9, 65, 61) & 585 \\
600 & 20 & (3, 8, 7) & 120 \\
624 & 18 & (13, 48, 43) & 624 \\
640 & 24 & (5, 8, 7) & 160 \\
840 & 28 & (24, 35, 31) & 840 \\
864 & 24 & (3, 8, 7) & 144 \\
880 & 22 & (16, 55, 49) & 880 \\
945 & 21 & (5, 21, 19) & 315 \\
945 & 27 & (7, 15, 13) & 315 \\
1000 & 30 & (5, 8, 7) & 200 \\
1056 & 16 & (11, 96, 91) & 1056 \\
1080 & 30 & (8, 15, 13) & 360 \\
1120 & 20 & (7, 40, 37) & 560 \\
1176 & 28 & (3, 8, 7) & 168 \\
1320 & 36 & (33, 40, 37) & 1320 \\
1344 & 36 & (16, 21, 19) & 672 \\
1360 & 24 & (17, 80, 73) & 1360 \\
1440 & 36 & (5, 8, 7) & 240\rlap{\enspace\yesmark}
\end{tabular}
\end{center}
\end{table}
The last line shown corresponds to the tiling in Fig.~\ref{figure:bryce1440}.  It
is not known if any of the other lines correspond to tilings.  We note that
while Theorem~\ref{theorem:abdividesX} has been formalized in Lean,
the correctness of the code that uses that theorem to compute Table~2 has not 
been formalized, so Table~2 is not a formal result.

Table~\ref{table:2} was computed as follows.  For each $N$ and each $M$
with $1 \le M < \sqrt N$, test whether $(9N-M^2)(N-M^2)$ is a square $r^2$
(it must be, by Lemma~\ref{lemma:square60}).  If so, by
Theorem~\ref{theorem:equilateral3} the tile $(a,b,c)$ is proportional to
$$(3N + M^2 - r,\ 3N + M^2 + r,\ 2(3N - M^2));$$
divide these by their greatest common divisor.  By the coloring equation,
the side of $ABC$ is $X = M(a+b+c)/3$.  Reject the line unless $ab$ divides
$X$ (Theorem~\ref{theorem:abdividesX}).  

The computation of the table is
instantaneous, and we could compute as many pages 
as we want to read.  

\begin{remark}
Laczkovich's method yields a necessary condition: the tiling must witness a relation
$jc = \ell a + mb$ with $j > 0$ and $\ell, m \ge 0$.  This does not help with the table,
because such a relation with $jc < X$ always exists.  Take $a, b, c$ to be integers with
no common factor (Theorem~3.3 of \cite{laczkovich2012}).  Then $a$ and $b$ are relatively
prime (as shown in the proof of Theorem~\ref{theorem:abdividesX}), and $a, b \ge 2$: if
$a = 1$, then $c^2 = b^2 \pm b + 1$ (the sign is $-$ for $\gamma = \pi/3$ and $+$ for
$2\pi/3$) lies strictly between $(b-1)^2$ and $(b+1)^2$ and differs from $b^2$, which is
impossible for $b \ge 2$, while $a = b = 1$ makes the tile equilateral; similarly if $b = 1$.
By a classical theorem of Sylvester \cite{sylvester1884}, every integer
$n \ge (a-1)(b-1)$ has the form $\ell a + mb$ with $\ell, m \ge 0$.  Let $j$ be
the least positive integer with $jc \ge (a-1)(b-1)$; then $jc = \ell a + mb$ for
some $\ell, m \ge 0$.  By the triangle inequality, $c \le a + b - 1$.  If $j = 1$,
then $jc = c \le a + b - 1 < ab$, since $(a-1)(b-1) \ge 1$.  If $j \ge 2$, then
$(j-1)c < (a-1)(b-1)$, so
$$jc < (a-1)(b-1) + c \le (a-1)(b-1) + a + b - 1 = ab.$$
In either case $jc < ab \le X$, by Theorem~\ref{theorem:abdividesX}.
\end{remark}

\noindent{\em Example}.  In Fig.~\ref{figure:bryce1440}, we have
$N = 1440$, and $(a,b,c) = (5,8,7)$, and $X = 240$.  Then we have
the relations
$5c = 7a$ and $5b = 8a$, $8c=7b$, and $8b = 3a + 7c$.
 Identify the segments in that tiling
that witness these relations.  Identify the graphs $\Gamma_c$ and $\Gamma_b$.
\medskip

\begin{lemma} \label{lemma:threeangles}
Suppose the equilateral triangle $ABC$ is tiled
by a tile with angles $(\alpha,\beta,\pi/3)$ with $\beta$ not
a rational multiple of $\pi$.  Then every vertex $V$ with angle sum $\pi$ has
three tiles with vertices at $V$, either all $\pi/3$ angles or one
each of $\alpha$, $\beta$, and $\pi/3$.
\end{lemma}

\begin{proof}  If there is no $\pi/3$ angle at $V$ then $\pi = p\alpha + q \beta$.
Since $\alpha + \beta = 2\pi/3$ we have $\pi = p(2\pi/3-\beta) + q \beta$. 
Then $\pi(1- 2p/3) = (q-p)\beta$.  Since $\beta$ is not a rational multiple of $\pi$,
we have $q=p$.  Then $\pi = p(\alpha + \beta) = 2p \pi/3$.  Hence $2p = 3$,
which is impossible.  Hence there is at least one $\gamma = \pi/3$ angle at $V$.
If there are exactly two $\gamma$ angles at $V$ then $\pi/3$ is a sum $p\alpha + q \beta$.
But one of $\alpha$ and $\beta$ is greater than $\pi/3$; hence $p=0$ or $q=0$.
Hence $\pi/3$ is a multiple of $\alpha$ or $\beta$, which is a contradiction.
More than three $\gamma$ angles adds to more than $\pi$.  
Therefore the number of $\gamma$ angles at $V$ is one or three.  If it is one,
then the remaining $2\pi/3$ can only be composed of one $\alpha$ and one $\beta$ 
angle, since if $p \alpha + q \beta = 2\pi/3$ then $\p(2\pi/3-\beta) + q \beta = 2\pi/3$,
so $2\pi/3(p-1) = (1-q)\beta$, which implies $q=p = 1$.
\end{proof}

\begin{lemma} \label{lemma:no24}
Suppose the equilateral $ABC$ is $N$-tiled by a
tile $(\alpha,\beta,\pi/3)$ with $\alpha$ not a rational multiple of $\pi$.
Then $N \ne 24$.
\end{lemma}

\begin{proof}  Suppose $N = 24$ and $ABC$ is $N$-tiled.  Let $M$ be the
coloring number of the tiling.  By Theorem~\ref{theorem:equilateral3}
and the coloring equation, $0 < M < \sqrt{24}$, so $1 \le M \le 4$.
By Lemma~\ref{lemma:square60}, $(9N-M^2)(N-M^2)$ is a square, and
only $M=4$ makes that expression a square, so $M=4$.
By Theorem~\ref{theorem:equilateral3}, part (ii), the tile is $(a,b,c) = (3,8,7)$,
and $ABC$ has side $X = 24$. Since $\alpha + \beta = 2\pi/3$ and $\alpha$ is
not a rational multiple of $\pi$, no sum of $\alpha$'s and $\beta$'s is $\pi/3$;
hence at each vertex of $ABC$ there is exactly one tile, with its $\gamma$
angle there.  In particular no tile has an $\alpha$ or $\beta$ angle at a
vertex of $ABC$.

Let Tile~1 be the tile at $A$.  Its edges at $A$ are its $a$ and $b$ edges, and
they lie along $AB$ and $AC$.  Renaming $B$ and $C$ if necessary, we may assume
its $b$ edge $AV$ lies on $AB$ and its $a$ edge $AX$ on $AC$.  Since the $b$ edge
is opposite the $\beta$ angle, Tile~1 has its $\alpha$ angle at $V$, its $\beta$
angle at $X$, and its $c$ edge is $VX$.  Since $8 < 24$ and $3 < 24$, $V$ lies
strictly inside $AB$ and $X$ strictly inside $AC$.

The angle sum at $V$ is $\pi$: no tile has $V$ in the interior of one of its
edges, since that tile would contribute $\pi$ to the angle sum at $V$, and
Tile~1 contributes $\alpha$ more.  By Lemma~\ref{lemma:threeangles}, there are
exactly three tiles with a vertex at $V$, with one $\alpha$, one $\beta$, and
one $\gamma$ angle there.  The same holds at $X$.

Consider the tiles on the other side of $VX$ from Tile~1 that have an edge
along $VX$.  Their edges cover $VX$, and none extends beyond $V$ or $X$, since
$V$ and $X$ are not in the interior of any tile edge.  So their lengths, each
$3$, $7$, or $8$, add up to $7$.  That can only be a single edge of length $7$.
Hence there is a single tile, Tile~2, whose $c$ edge is $VX$.  Its $\gamma$
angle is opposite $VX$, so at $V$ it has $\alpha$ or $\beta$.  It is not
$\alpha$, since Tile~1 already has the only $\alpha$ at $V$.  Hence Tile~2 has
its $\beta$ angle at $V$.

We claim that no tile with an edge on $AB$ has its $\beta$ angle at $V$.  Such
a tile would be Tile~2, since there is only one $\beta$ at $V$.  Then Tile~2
would have an edge $VY$ on $AB$, besides its edge $VX$.  The point $Y$ cannot
lie between $A$ and $V$, since then the edges $VY$ and $AV$ of Tiles~2 and~1,
both on $AB$, would overlap.  So $V$ lies between $A$ and $Y$, and the angle
of Tile~2 at $V$ is $\angle XVY = \pi - \angle AVX = \pi - \alpha$.  But that
angle is $\beta$, so $\alpha + \beta = \pi$, contradicting
$\alpha + \beta = 2\pi/3$.  That proves the claim.

Now let $n$ tiles have edges on $AB$; of these, let $i$ have their $a$ edges
there, $j$ their $b$ edges, and $k$ their $c$ edges.  Then $n = i + j + k$ and
$$ 3i + 8j + 7k = 24.$$
These $n$ edges have $n-1$ endpoints strictly inside $AB$.  A tile with its
$a$ or $c$ edge on $AB$ has a $\beta$ angle at an endpoint of that edge.  That
endpoint is not $A$ or $B$, and by the claim it is not $V$.  By
Lemma~\ref{lemma:threeangles}, no vertex inside $AB$ has two $\beta$ angles.
Hence $i + k \le n - 2$, that is, $j \ge 2$.  Similarly, a tile with its $b$
or $c$ edge on $AB$ has an $\alpha$ angle at an endpoint of that edge, not at
$A$ or $B$, and no two at the same vertex; hence $j + k \le n - 1$, that is,
$i \ge 1$.  Then $8j \le 24 - 3i \le 21$, so $j = 2$ and $3i + 7k = 8$.  But
$3i + 7k = 8$ has no solution in nonnegative integers.  This contradiction
completes the proof.
\end{proof}

We turn now to necessary conditions in case $\gamma = 2\pi/3$.  Then
there is no coloring number, but we can still find the possible $(a,b,c)$
when given $N$, since $ab$ must divide $N$.  For efficiency we use 
the equation $c^2 = a^2 + b^2 + ab$,  whose solutions are parametrized by 
$$ (a,b,c) = (s^2-t^2,2st+t^2,s^2+st+t^2),  $$
where $s > t > 0$, $\gcd(s,t) = 1$, and $s \not\equiv t \pmod 3$, then
reordered so $a < b$.  For example, $(s,t) = (2,1)$ gives $(3,5,7)$, 
$(3,1)$ gives $(7,8,13)$ after swapping, and $(3,2)$ gives $(5,16,19)$.
One can find this parametrization, for example, in \cite{cohen2007}.

We reject lines in which $ab$ does not divide $N$.  We also 
rejected $N=15$, which would require $X = 15$ and $(a,b,c) = (3,5,7)$, on the 
grounds that a side of length 15 cannot support two tiles with $c$ edges 
of $7$,  contrary to Lemma~8 of \cite{beeson-triangletiling3}.  (That is 
the only line eliminated from the table by that lemma.)

\begin{table} [ht]
\caption{Tilings with $\gamma = 2\pi/3$ and $N \le 1500$ not ruled out by \\
Corollary~\ref{corollary:abdividesN} and the two-$c$-edges lemma.}
\label{table:3}
\begin{center}
\begin{tabular}{rrrr}
$N$ &  $m$ & $(a,b,c)$ & side of $ABC$ \\
\hline
56 & 1 & (7, 8, 13) & 56 \\
60 & 2 & (3, 5, 7) & 30 \\
80 & 1 & (5, 16, 19) & 80 \\
135 & 3 & (3, 5, 7) & 45 \\
224 & 2 & (7, 8, 13) & 112 \\
231 & 1 & (7, 33, 37) & 231 \\
240 & 4 & (3, 5, 7) & 60 \\
264 & 1 & (11, 24, 31) & 264 \\
320 & 2 & (5, 16, 19) & 160 \\
375 & 5 & (3, 5, 7) & 75 \\
455 & 1 & (13, 35, 43) & 455 \\
504 & 1 & (9, 56, 61) & 504 \\
504 & 3 & (7, 8, 13) & 168 \\
540 & 6 & (3, 5, 7) & 90 \\
624 & 1 & (16, 39, 49) & 624 \\
720 & 3 & (5, 16, 19) & 240 \\
735 & 7 & (3, 5, 7) & 105 \\
896 & 4 & (7, 8, 13) & 224 \\
924 & 2 & (7, 33, 37) & 462 \\
935 & 1 & (11, 85, 91) & 935 \\
960 & 8 & (3, 5, 7) & 120 \\
1056 & 2 & (11, 24, 31) & 528 \\
1071 & 1 & (17, 63, 73) & 1071 \\
1215 & 9 & (3, 5, 7) & 135 \\
1280 & 4 & (5, 16, 19) & 320 \\
1400 & 5 & (7, 8, 13) & 280 \\
1440 & 1 & (32, 45, 67) & 1440 \\
1500 & 10 & (3, 5, 7) & 150
\end{tabular}
\end{center}
\end{table}

\section{Appendix: Tilings found by Bryce Herdt}
In 2024, Bryce Herdt found new tilings of equilateral triangles 
using tiles $(3,8,7)$ and $(5,8,7)$ (for which $\gamma = \pi/3$),
and $(3,5,7)$ (for which $\gamma = 2\pi/3$).  These tilings dramatically lowered the $N$
for the ``smallest known tiling'' of equilateral triangles. 
Herdt found these tilings by first finding new dissections of an equilateral triangle into
similar triangles, and then refining those to tilings by congruent tiles.  

\begin{figure}[ht]
\caption{Decomposing a parallelogram with top $pc + qb$ and side $bc$}
\label{figure:bryce1}
\begin{center}
\includegraphics[width=0.8\textwidth]{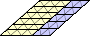}
\end{center}
\end{figure}

Herdt's key 
innovation is the realization that a parallelogram can often be broken into two parallelograms,
which can then be tiled with tiles in different orientations, dramatically reducing the 
number of tiles required. 
 Fig.~\ref{figure:bryce1}
illustrates the technique. 

\begin{figure}[ht]
\caption{Herdt's 1944-tiling with tile $(3,8,7)$}
\label{figure:bryce1944}
\begin{center}
\includegraphics[width=\textwidth]{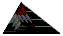}
\end{center}
\end{figure}

\begin{figure}[ht]
\caption{Herdt's 1440-tiling with tile $(5,8,7)$}
\label{figure:bryce1440}
\begin{center}
\includegraphics[width=\textwidth]{Figures/tiling1440}
\end{center}
\end{figure}


\clearpage
%

%

\appendix
\section{Formal verification in Lean}\label{section:lean}

The results of this paper have been formalized in the Lean~4 proof assistant
(version~4.30.0), using the Mathlib library.  The Lean code is at
\url{https://github.com/mjbeeson/TriangleTilingEquilateralLean}, release v1.0.
Each theorem, lemma and corollary listed in Table~\ref{table:lean} has a Lean
counterpart, named after its label, which states the result as printed here and
whose proof follows the printed proof.  The proofs are checked by the Lean kernel
using only the standard axioms ({\tt propext}, \goodbreak\noindent {\tt Classical.choice},
{\tt Quot.sound}), with no {\tt sorry} and no {\tt native\_decide}.
All Lean names are in the namespace {\tt TriangleTiling.Equilateral}, except
\goodbreak\noindent
{\tt notallsupported}, which is in {\tt TriangleTiling}.

Tilings are formalized as in V.~Bonfioli's Lean library \cite{bonfioli-lean}:
$N$ triangles congruent to a given tile (reflections allowed), with disjoint
interiors, whose union is $ABC$.  Two tiles {\em share a segment} if an edge of
each contains a common segment of positive length.  A path from $A$ that crosses
tile boundaries away from vertices is formalized as a chain of tiles, each sharing
a segment with the next; the number of boundaries crossed is the length of the
chain.  ``The tile at $A$'' is a tile with a vertex at $A$ and an edge on $AB$,
and the coloring number $M$ is the sum of a proper coloring of the tiles by
$\pm 1$ with that tile colored $+1$.

\smallskip
\noindent{\em What is not formalized.}  Tables~\ref{table:2} and~\ref{table:3},
which are computations; every condition used to compute them is formalized.  The
general form of the coloring theorem in the Remark after
Theorem~\ref{theorem:coloring}, and the search procedure in the Remark after
Theorem~\ref{theorem:decidable}.

\smallskip
\noindent{\em Results taken from Bonfioli's library.}  Laczkovich's results used
here are taken from \cite{bonfioli-lean}: Theorem~3.3 of \cite{laczkovich2012}
(the sides of the tile are commensurable), Lemma~3.2 of \cite{laczkovich2012}
(the tiling is regular), and the definitions and basic properties of the graph
$\Gamma_c$, from which Lemmas~\ref{lemma:propagation2}
and~\ref{lemma:noboundarylink} are proved.  So is the fact, used in the proof of
Theorem~\ref{theorem:coloring}, that every tile can be reached from the tile at $A$
by a chain of tiles sharing segments.  The results marked $\dagger$ in
Table~\ref{table:lean} depend on these and are in a separate part of the Lean
code, which takes several hours to build; the README of the repository explains
how.

\begin{table}[ht]
\caption{Results of this paper and their Lean names ($\dagger$: in the separately
built part).}
\label{table:lean}
\begin{center}
\begin{tabular}{ll}
Result & Lean name \\
\hline
Theorem~\ref{theorem:coloring} & {\tt coloring}$^\dagger$ \\
Theorem~\ref{theorem:equilateral3} & {\tt equilateral3}$^\dagger$ \\
Lemma~\ref{lemma:ab} & {\tt ab}$^\dagger$ \\
Lemma~\ref{lemma:square60} & {\tt square60}$^\dagger$ \\
Theorem~\ref{theorem:notprime60} & {\tt notprime60}$^\dagger$ \\
Theorem~\ref{theorem:decidable} & {\tt decidable}$^\dagger$ \\
Theorem~\ref{theorem:comparison} & {\tt comparison}$^\dagger$ \\
Lemma~\ref{lemma:propagation2} & {\tt propagation2}$^\dagger$ \\
Lemma~\ref{lemma:noboundarylink} & {\tt noboundarylink}$^\dagger$ \\
Lemma~\ref{lemma:notallsupported} & {\tt notallsupported} \\
Theorem~\ref{theorem:notprime120} & {\tt notprime120}$^\dagger$ \\
Theorem~\ref{theorem:abdividesX} & {\tt abdividesX} \\
Corollary~\ref{corollary:abdividesN} & {\tt abdividesN} \\
Lemma~\ref{lemma:threeangles} & {\tt threeangles} \\
Lemma~\ref{lemma:no24} & {\tt no24}$^\dagger$ \\
\end{tabular}
\end{center}
\end{table}

\vskip 20pt
\noindent {\bf Acknowledgment.} The original version of this paper was
written without AI assistance.  In preparing this revision, the author
used Claude (Anthropic), an AI system, in three ways.  First, to check
the proofs: it found errors and gaps, among them a sign error in the
original proof of Theorem~\ref{theorem:notprime60}, an off-by-one error in
the computation of Table~\ref{table:2}, an unproved claim in the proof of
Theorem~\ref{theorem:coloring} (that the coloring does not depend on the
path), unjustified steps in the original proof of
Theorem~\ref{theorem:equilateral3}, and missing cases in the proof of
Theorem~\ref{theorem:notprime120}; and it proposed corrections. Second, it proposed simpler proofs of
Theorem~\ref{theorem:equilateral3} and Lemma~\ref{lemma:no24}, and the
argument for path independence in Theorem~\ref{theorem:coloring}.
 Third,
it wrote the Lean formalization described in Appendix~\ref{section:lean},
whose proofs are checked by the Lean kernel.  The author checked all
proposed changes to the mathematics, rejected or modified some, checked
that the Lean statements say what the theorems of the paper say, and is
responsible for the final text.  The main results are unchanged.


\begin{thebibliography}{1}

\bibitem{beeson-noseven}
Michael Beeson.
\newblock {\em No triangle can be cut into seven congruent triangles}.
\newblock 2018.
\newblock \url{https://arxiv.org/abs/1811.09723}.

\bibitem{beeson-triangletiling3}
Michael Beeson.
\newblock Triangle tiling: The case $3\alpha + 2\beta = \pi$, 2019.
\newblock \href {http://arxiv.org/abs/1206.2229} {\path{arXiv:1206.2229}}.

\bibitem{bonfioli-lean}{\raggedright V. Bonfioli,
\newblock Formalizing triangle tiling in Lean, version 1.0, 22 October 2026.
\newblock Lean~4 formalization, GitHub repository
\texttt{ElVec1o/formalizing-triangle-tiling-in-lean}.\par}


\bibitem{cohen2007}
H.~Cohen.
\newblock {\em Number Theory, Volume I: Tools and Diophantine Equations}.
\newblock Graduate Texts in Mathematics 239, Springer, New York, 2007.

\bibitem{harries2026no19tiling}
Jan Philipp Harries.
\newblock{\em  No triangle can be cut into nineteen congruent triangles: The prime case of \Erdos problem~634}.
\newblock \url{https://github.com/jphme/math-problems/tree/no19tiling-v2/no19tiling}.

\bibitem{harries2026}
Jan Philipp Harries.
\newblock{\em  New constructions, obstructions, and multiplier structure for \Erdos problem~634.}
\newblock \url{https://github.com/jphme/math-problems/tree/main/progress634}, v.~5.


\bibitem{laczkovich1995}
M.~Laczkovich.
\newblock {Tilings of triangles}.
\newblock {\em Discrete Mathematics}, 140:79--94, 1995.


\bibitem{laczkovich2012}
Mikl{\'o}s Laczkovich.
\newblock Tilings of convex polygons with congruent triangles.
\newblock {\em Discrete and Computational Geometry}, 38:330--372, 2012.


\bibitem{snover1991}
Stephen~L. Snover, Charles Waiveris, and John~K. Williams.
\newblock {Rep-tiling for triangles}.
\newblock {\em Discrete Mathematics}, 91:193--200, 1991.


\bibitem{soifer}
Alexander Soifer.
\newblock {\em How Does One Cut a Triangle?}
\newblock Springer, 2009.


\bibitem{Sage}
The Sage Developers.
\newblock {\em {S}age{M}ath, the {S}age {M}athematics {S}oftware {S}ystem
  ({V}ersion 10.9)}, 2026.
\newblock \url{https://www.sagemath.org}.

\bibitem{sylvester1884}
J. J. Sylvester, J. J,
\newblock Mathematical questions with their solutions.
\newblock {\em Educational Times} 41:21, 1884. 

\bibitem{zhang2026}
Yan~X Zhang.
\newblock Tiling triangles with $2\pi/3$ angles, 2026.
\newblock \href {http://arxiv.org/abs/2512.22696v4} {\path{arXiv:2512.22696}}.



\end{thebibliography}
\bibliographystyle{plain-annote}

\end{document}